\documentclass{amsart}

\newtheorem{theorem}{Theorem}[section]

\theoremstyle{definition}
\newtheorem{definition}[theorem]{Definition}

\theoremstyle{remark}

\numberwithin{equation}{section}

\usepackage{hyperref}

\begin{document}

\title{Elementary Evaluation of the Zeta and Related Functions: An Approach From a New Perspective}

\author{Armen Bagdasaryan}
\address{Institution of the Russian Academy of Sciences, Trapeznikov Institute for Control Sciences of RAS, 
65 Profsoyuznaya, 117997, Moscow, Russia}
\email{abagdasari@hotmail.com}



 



\begin{abstract}
The zeta and related functions are explicitly known at either even values as for the $\zeta(s)$, $\lambda(s)$ and $\eta(s)$ functions
or at odd function values as for the Dirichlet $\beta(s)$ function. In this paper, the values at integers are obtained for the zeta and
related functions by applying new analytical tool that naturally arises within a new theoretical setting. This allows a simple
and elementary derivation of the function values at integer points from quite a different point of view. It is based on the
hypothesis that negative numbers might be beyond infinity, advanced by Wallis and Euler. I will first define the basic concepts
of our theory: (1) a new ordering relation on the set of integer numbers; (2) a new class of real functions, called regular, and
the definition of sum that extends the classical definition to the case when the upper limit of summation is less than the lower;
(3) a set of conditions imposed on regular functions that define a new regular method for summation of infinite series.
\end{abstract}

\maketitle

\section{The Riemann's Zeta and Related Functions and their Evaluations}
The evaluations of series associated with the Riemann's zeta and related
functions have a long history which can be traced back to Christian Goldbach
(1690--1764) and Leonhard Euler (1707--1783) \cite{monthly,kline,vara}.
Since then, many different techniques to evaluate various families of series involving the
zeta and related functions have been developed \cite{apostol}--\cite{borwein}.

The main object of this paper is to show how unexpectedly simply and
nicely 
the Riemann's zeta function and some related functions can be
evaluated from quite a different point of view by starting with a single basic concept.

The Riemann zeta function is defined by the Dirichlet series

$$
\zeta(s)=\sum_{n=1}^{\infty}\frac{1}{n^s}
$$

convergent for $\Re(s)>1$, which is guaranteed by the integral test. The function has an analytic continuation to the
whole complex plane with a simple pole at $s=1$ with residue $1$.
The values of $\zeta(s)$ at integer arguments have significant arithmetic meanings, and are of much interest not only
in number theory, but also in physics and statistics.

One of the most celebrated formulas, discovered by Euler in 1734, is the formula for even zeta values

\begin{equation} \label{zetaeven}
\zeta(2k)=(-1)^{k}\;\frac{2^{2k-1}\,B_{2k}}{(2k)!}\;\pi^{2k}
\end{equation}

in which $k$ is a positive integer and $B_{k}$ are the Bernoulli numbers. No analogous closed form representation of $\zeta(s)$
at odd integers is known. However, a number of series representations, as well as several integral representations, for the values of
$\zeta(2k+1)$ have been derived \cite{srivajour1,sriva}.

There are functions that are closely related to $\zeta(s)$. Namely, the Dirichlet eta $\eta(\cdot)$ and lambda $\lambda(\cdot)$ functions
that are given by

\begin{align} \label{eta}
\eta(s)=\sum_{n=1}^{\infty}\frac{(-1)^{n-1}}{n^s}=\frac{1}{\Gamma(s)}\int_{0}^{\infty}\frac{t^{s-1}}{e^{t}+1}dt, \quad s>0 \nonumber \\ 
 \\
\lambda(s)=\sum_{n=0}^{\infty}\frac{1}{(2n+1)^s}=\frac{1}{\Gamma(s)}\int_{0}^{\infty}\frac{t^{s-1}}{e^{t}-e^{-t}}dt, \quad s>0 \nonumber
\end{align}

where $\Gamma(s)$ is the Euler's gamma function.
These functions are related to $\zeta(s)$ by

\begin{equation} \label{relation}
\eta(s)=(1-2^{1-s})\zeta(s) \quad \textup{and} \quad \lambda(s)=(1-2^{-s})\zeta(s)
\end{equation}

They satisfy the simple identity $\zeta(s)+\eta(s)=2\lambda(s)$. As a consequence of the relation to $\zeta(s)$
via (\ref{relation}), the explicit expressions for both of $\eta(2k)$ and $\lambda(2k)$ exist.

Another function related to $\zeta(s)$ is the Dirichlet beta $\beta(\cdot)$ function. This function is defined as

\begin{equation} \label{beta}
\beta(s)=\sum_{n=0}^{\infty}\frac{(-1)^n}{(2n+1)^s}=\frac{1}{\Gamma(s)}\int_{0}^{\infty}\frac{t^{s-1}}{e^{t}+e^{-t}}dt, \quad s>0
\end{equation}

where $\beta(2)=G$ is the Catalan's constant, a series that converges for all $s>0$, accordingly to the Leibnitz's test for
alternating series. It can be observed from (\ref{eta}) that $\beta(s)$ is the alternating version of $\lambda(s)$. 
However, $\beta(s)$ cannot be directly related to $\zeta(s)$, but it is related to $\eta(s)$ in that only the odd terms are summed.

The function $\beta(s)$ can be explicitly evaluated at positive odd integer points of $s$, as follows

$$
\beta(2k+1)=(-1)^{k}\;\frac{E_{2k}}{2^{2k+2}\,(2k)!}\;\pi^{2k+1}
$$

where $E_{k}$ are the Euler numbers. 

As for the mysterious situation around $\zeta(2k+1)$,  the same scenario we have 
for even values of $\beta(s)$, that is no finite closed form expression for $\beta(2k)$ is known. In this connection, 
$\beta(2k)$ can be considered as a ``reverse'' counterpart to the mystery of $\zeta(2k+1)$.
The analytic continuation extends the function $\beta(s)$ to all points in the complex plane. This can be done by the functional equation

$$
\beta(1-s)=\left(\frac{2}{\pi}\right)^{s}\,\sin\left(\frac{\pi s}{2}\right)\,\Gamma(s)\,\beta(s)
$$

For negative integer values of the argument, this yields $\beta(-k)=E_{k}/2$. Unlike the Riemann's $\zeta(s)$, the Dirichlet's 
$\beta(s)$ has no singularities.

The values of $\zeta(s)$ at negative integers is known to be expressed via the Bernoulli numbers, $\zeta(-k)=-B_{k+1}/(k+1)$. 
By the relation (\ref{relation}), the finite formulas for $\eta(-k)$ and $\lambda(-k)$ are readily derived. 

In deriving above values, basically at negative integers, most of the methods employ complex-analytical tools, the techniques of
analytic continuation and contour integration are often used. In contrast to these methods, in this work I present a new elementary derivation 
of these values, based on the new theoretical standpoints.

\section{Basic Definitions and Axiomatics for Summation of Series}

A new ordering relation on the set of all integer numbers $\mathbb{Z}$ is introduced as follows \cite{physnuclei,jmr,varbag}. 

\begin{definition}
Let $a, b \in \mathbb{Z}$. Then $a$ precedes $b$, $a\prec b$, if and only if \, $\frac{-1}{a}<\frac{-1}{b}$; $a\prec b \Leftrightarrow \frac{-1}{a}<\frac{-1}{b}$ \footnote{by convention $0^{-1}=\infty$}.
\end{definition}

Hence, we get that $\mathbb{Z}=[0,1,2,3,...,-3,-2-1]$. The new ordering satisfies the axioms of transitivity and connectedness, 
so it defines a strict total (linear) order on $\mathbb{Z}$. 

Let $f(x)$ be a function of real variable defined on $\mathbb{Z}$. 

\begin{definition} \label{regdef}
The function $f(x)$, $x\in \mathbb{Z}$, is called \textit{regular} if there exists an elementary\footnote{determined by formulas, constructed by a finite number of algebraic operations and constant functions and algebraic, trigonometric, exponential, and logarithmic functions, and their inverses through repeated combinations and compositions}
function $F(x)$ such that $F(z+1)-F(z)=f(z), \; \forall z\in \mathbb{Z}$. 
\end{definition}

From Def. \ref{regdef} follows that $F(b+1)-F(a)=\sum_{u=a}^{b}f(u)$.
Let $\mathbb{Z}_{a, b}$ denotes a part of $\mathbb{Z}$ such that $\mathbb{Z}_{a, b}=[a,b]$ if $(a\preceq b)$ and 
$\mathbb{Z}_{a, b}=\mathbb{Z}\setminus (b, a)$ if $(a\succ b)$, where $\mathbb{Z}\setminus (b, a)=[a, -1]\cup[0, b]$.
Then we introduce the following fundamental

\begin{definition}[sums of regular functions]
$
\sum_{u=a}^b{f(u)}=\sum_{u\:\in\: \mathbb{Z}_{a, b}}{f(u)}, \; \forall a, b\in \mathbb{Z} 
$
\end{definition}
which extends the usual sum.
 
We impose some quite natural conditions on regular functions. Namely,

\begin{enumerate}\baselineskip 16pt
\item[\bf A1.] If \, $S_n=\sum_{u=a}^n{f(u)} \quad \forall n$, \quad then \, $\lim_{n\rightarrow \infty}S_n=\sum_{u=a}^\infty {f(u)}$. 
\item[\bf A2.] If \, $S_n=\sum_{u=1}^{n/2}{f(u)} \quad \forall n$, \quad then \, $\lim_{n\rightarrow \infty}S_n=\sum_{u=1}^\infty {f(u)}$.
\item[\bf A3.] If \, $\sum_{u=a}^\infty {f(u)}=S$, \quad then \, $\sum_{u=a}^\infty {\lambda f(u)}=\lambda S, \; \lambda\in \mathbb{R}$.
\item[\bf A4.]  If \, $\sum_{u=a}^\infty {f_1(u)}=S_1$ \, and \, $\sum_{u=a}^\infty {f_2(u)}=S_2$, \quad then \, $\sum_{u=a}^\infty {\left(f_1(u)+f_2(u)\right)}= S_1+S_2$.
\item[\bf A5.] If \, $G=[a_1, b_1]\cup[a_2, b_2]$, $[a_1, b_1]\cap[a_2, b_2]=\emptyset$, \quad then \, $\sum_{u\in G}{f(u)}=\sum_{u=a_1}^{b_1}{f(u)}+\sum_{u=a_2}^{b_2}{f(u)}$. 
\end{enumerate}
 
This axiomatic system defines the method of summation of infinite series, which is regular due to \textbf{A1}.

\begin{definition}\label{quasi}
The function $f(x)$, $x\in\mathbb{R}$, is called \textit{quasi-even} if the condition $f(-x)=f(x-a)$, $a\in\mathbb{Z}$, holds.
\end{definition}

One of the formulas that we have at the core of our theory is given by the following theorem.

\begin{theorem}[summation formula]
\label{theoremsum}
Let $f(x)$ be a regular quasi-even function that satisfies the condition of definition \ref{quasi} with $a=\epsilon t$, 
where $\epsilon=\pm 1$ and $t$ is a fixed natural number, and let \, $\delta=2^{-1}(1-\epsilon)$. Then, it holds that
\begin{equation}
\sum_{u=1}^{\infty}f(u)\;=\;\frac{1}{2}\,\epsilon\;\sum_{u=\delta}^{t-1+\delta}\left(\lim_{n\to\infty}f\left(n-\epsilon u\right)-f\left(-\epsilon u\right)\right)-\frac{1}{2}f(0).
\label{summation}
\end{equation}
\end{theorem}

Formally, for $\epsilon=0$, formula \ref{summation} reduces to \; $\sum_{u=1}^{\infty}f(u)=-f(0)/2$ and coincides with the 
statement of theorem 6 in \cite{physnuclei}, 
in which $f(x)$ is supposed to be even. 
So, taking \, $F(n)=B_{2k}(n-1)$, where $B_{2k}(x)$ is the Bernoulli polynomial, we immediately deduce \; $\sum_{u=1}^{\infty}u^{2k}=0 \quad \forall k$, that is, the 
trivial zeroes of $\zeta(s)$, the values of $\zeta(-2k)$.

\section{Application to Evaluating of the Zeta Functions}

\subsection{Values at non-positive integers.}
The values of $\zeta(s)$ and $\eta(s)$ have recently been obtained in \cite{physnuclei} by means of the techniques of our new setting.
The explicit formulas for the sum of infinite arithmetic and alternating arithmetic progression, as well as some other original results, were used there
to derive the values of $\zeta(-k)$ and $\eta(-k)$. Using the relation (\ref{relation}), the values of $\lambda(-k)$ are derived.
By means of similar techniques, the values of $\beta(s)$ at negative integers can be obtained. Using formulas (8) and (9) from \cite{physnuclei} 
along with some other formulas, we get 
$$
\beta(1-k)=\sum_{u=1}^{\infty}(-1)^{u-1}(2u-1)^{k-1}=-\frac{1}{2k}\,\sum_{u=1}^{k}(-1)^{u}\,2^{u}\left(2^{u}-1\right)\binom{k}{u}B_{u} \qquad \forall \, k\in\mathbb{N}.
$$
In the same manner, we can obtain the values of $\lambda(1-k)$.
\newline

\subsection{Values at positive integers.}
Let's take the Taylor expansion of the sine function, 
$$
\sin x = \sum_{u=0}^{\infty}\frac{(-1)^{u}}{(2u+1)!}x^{2u+1} = x - \frac{x^3}{3!} + \frac{x^5}{5!} - ... \qquad |x|<\infty
$$
Since convergent series can be multiplied by constants and be term-by-term added (subtracted), the series
\begin{eqnarray}
\sum_{u=1}^{\infty}(-1)^{u-1}\frac{\sin ux}{u^{2k+1}}= x\sum_{u=1}^{n}\frac{(-1)^{u-1}}{u^{2k}}-\frac{x^3}{3!}\sum_{u=1}^{n}\frac{(-1)^{u-1}}{u^{2(k-1)}}+...+  \\
 (-1)^{k}\frac{x^{2k+1}}{(2k+1)!}\sum_{u=1}^{n}(-1)^{u-1} +  
(-1)^{k+1}\frac{x^{2k+3}}{(2k+3)!}\sum_{u=1}^{n}(-1)^{u-1}u^2 +... \nonumber
\end{eqnarray}
converges for all $n,k$ and $x$, $-\pi\leq x\leq \pi$. Using Rule III and formula (3.3) \cite{varbag}, we get
\begin{equation}\label{eqsine}
\sum_{u=1}^{\infty}(-1)^{u-1}\frac{\sin ux}{u^{2k+1}} = \sum_{\nu=0}^{k}(-1)^{\nu}\frac{x^{2\nu+1}}{(2\nu+1)!}\sum_{u=1}^{\infty}\frac{(-1)^{u-1}}{u^{2(k-\nu)}},
\qquad -\pi\leq x\leq \pi
\end{equation}
Putting in (\ref{eqsine}) $x=\pi$, we have
$$
0 = \sum_{\nu=0}^{k}(-1)^{\nu}\frac{\pi^{2\nu}}{(2\nu+1)!}\sum_{u=1}^{\infty}\frac{(-1)^{u-1}}{u^{2(k-\nu)}}
$$
Letting $k=1,2,3...$, we find 
$$
\sum_{u=1}^{\infty}\frac{(-1)^{u-1}}{u^2}=\frac{\pi^2}{12}, \qquad \sum_{u=1}^{\infty}\frac{(-1)^{u-1}}{u^4}=\frac{7\pi^4}{720},...
$$
and applying the induction over $k$, and using (4.3) \cite{varbag}, we obtain
\begin{equation}\label{eta:even}
\eta(2k) = \sum_{u=1}^{\infty}\frac{(-1)^{u-1}}{u^{2k}} = (-1)^{k-1}\frac{\left(2^{2k-1}-1\right)\pi^{2k}}{(2k)!}B_{2k}
\end{equation}

By the relations (\ref{relation}) the values of $\zeta(2k)$ and $\lambda(2k)$ are easily derived from (\ref{eta:even}).
Analogously, taking the Taylor expansion of the cos function, we have another way to derive the values of $\lambda(2k)$, and
then using (\ref{relation}), we get $\zeta(2k)$ and $\eta(2k)$.

To evaluate $\beta(s)$ at positive odd integers, we substitute the formula (\ref{eta:even}) into (\ref{eqsine}), and putting $\nu=k-u$, we have

\begin{align}\label{eqsine1}
\sum_{u=1}^{\infty}(-1)^{u-1}\frac{\sin ux}{u^{2k+1}} = \frac{(-1)^{k+1}}{(2k+1)!}\sum_{u=0}^{k}\left(2^{2u-1}-1\right)\binom{2k+1}{2u}
x^{2(k-u)+1}\pi^{2u}B_{2u}, \\
 -\pi\leq x \leq \pi   \nonumber
\end{align}

Then, replacing $x$ by $\pi-x$ in (\ref{eqsine1}), we get

\begin{align}\label{eqsinefinal}
\sum_{u=1}^{\infty}\frac{\sin ux}{u^{2k+1}} = \frac{(-1)^{k+1}}{(2k+1)!}\sum_{u=0}^{k}\left(2^{2u-1}-1\right)\binom{2k+1}{2u}
(\pi-x)^{2(k-u)+1}\pi^{2u}B_{2u},  \\
0\leq x \leq 2\pi  \nonumber
\end{align}

Now taking $x=\pi/2$ in (\ref{eqsinefinal}), we obtain

\begin{align}
\beta(2k+1)=\sum_{u=1}^{\infty}\frac{(-1)^{u-1}}{(2u-1)^{2k+1}} =  \frac{(-1)^{k+1}\pi^{2k+1}}{2^{2k+1}(2k+1)!}\sum_{u=0}^{k}2^{2u}\left(2^{2u-1}-1\right) \times \\ \times
\binom{2k+1}{2u}B_{2u} \nonumber
\end{align}

The same reasoning allows us to obtain some other results related to summation of series.

\section{Concluding Remarks}

In this work we have presented the basic concepts of the new theory and
its application to the evaluation of the zeta and related functions. 
More detailed presentation of the results along with the complete proofs of propositions will be published elsewhere.
In evaluating these functions, we have not employed the analytic continuation or
contour integration techniques, unlike most of the works that deal with the problem.
It would be interesting to develop some similar results and formulas for some other series or classes of series 
associated with the zeta and related functions.
The theoretical basics we outlined provides 
not only a new method to discover new results, but generally opens a new research direction; 
to the author's opinion, studying the new number line in physical context is of particular interest.

\bibliographystyle{amsplain}

\end{document}